\documentclass[11pt] {article} 

\pdfoutput = 1

\usepackage[backref]{hyperref}
\usepackage{amsmath}
\usepackage{amssymb}    
\usepackage{pdfsync}

\newtheorem{theorem}{Theorem}[section]
\newtheorem{lemma}{Lemma}[section]

\newcommand{\eqnsection}{
   \renewcommand{\theequation}{\thesection.\arabic{equation}}
   \makeatletter
   \csname @addtoreset\endcsname{equation}{section} 
   \makeatother}
    
\def \be{\begin{equation}}
\def \ee{\end{equation}}
\def \bt{\begin{theorem}} 
\def \et{\end{theorem}} 
\def \bl{\begin{lemma}}     
\def \el{\end{lemma}}
\def \bea{\begin{eqnarray}}
\def \eea{\end{eqnarray}}
\def \bas{\begin{eqnarray*}}
\def \eas{\end{eqnarray*}}

\def \de{\delta} 
\def \De{\Delta} 
\def \ep{\epsilon}

\def \la{\lambda} 

\def \ka{\kappa}
\def \om{\omega}
\def \Om{\Omega}

\def \si{\sigma}

\def \th{\theta}

\def \ze{\zeta}

\def \ff{\infty}
\def \wh{\widehat}
\def \wt{\widetilde}

\def \cd{\,\cdot\,}

\def \AA{{\cal A}}
\def \BB{{\cal B}}

\def \FF{{\cal F}}

\def \MM{{\cal M}}

\def \RR{{\cal R}}

\def \({\left(}
\def \){\right)}

\def \nn{\nonumber}

\def \bc{\begin{center} }
\def \ec{\end{center} }
\def \bs{\begin{slide} }
\def \es{\end{slide} }

\def\square{{\vcenter{\vbox{\hrule height.3pt
        \hbox{\vrule width.3pt height5pt \kern5pt
           \vrule width.3pt}
        \hrule height.3pt}}}}
\def\qed{{\hfill $\square$ \bigskip}}

 \def \Rev({\mbox{\rm Rev}(}

\eqnsection

\begin{document}
\title{Perturbation of the loop measure}

\author{Yves Le Jan \hskip20pt  Jay Rosen}

\maketitle

 \begin{abstract}
The loop measure is associated with a Markov generator. We compute the variation of the loop measure induced by an infinitesimal variation of the generator affecting the killing rates or the jumping rates.   
 \end{abstract}

\section{Introduction}
Professor Fukushima's contribution to probabilistic potential theory was  seminal. His book on Dirichlet spaces \cite{Fuk} was the first complete exposition in which the functional analytic, potential theoretic, and probabilistic aspects of the theory were considered jointly, as different aspects of the same mathematical object. In particular, transformations of the energy form, such as restriction to an open set, trace on a closed set, change of the equilibrium (or killing) measure, change of reference measure, etc... have probabilistic counterparts. In the second chapter, the possibility of superposing different Dirichlet forms was recognized. The present paper follows this line of research. But we will focus on Markov loops and bridges rather than Markov paths. We will present them in the next section. Let us stress however that the existence of loop or bridge measures seems to require more than the  assumption of a  Markov process defined up to polar sets, which is the basic assumption for a  Markov process associated with a Dirichlet form. The existence of a Green function seems to be required.
Our purpose is to compute the variation of the loop measure induced by an infinitesimal variation of the generator. This variation may in particular affect the killing rates or the jumping rates. In the case of symmetric continuous time Markov chains, the question has been addressed in chapter 6 of \cite{Le Jan1}. The results are formally close to formulas used in conformal field theory for operator insertions (see for example \cite{Gaw}). We try here to extend them to a more general situation, to show in particular that the bridge measures can be derived from the loop measure.

\section{Background on loop measures}

 We begin by introducing loop measures for Borel right processes (such as  Feller processes) on a rather general state space
    $S$, which we assume to  be locally compact   with a countable base. 
  Let     $X\!=\!
(\Om,  \FF_{t}, X_t,\th_{t},P^x
)$ be a transient Borel right process \cite{S} with cadlag paths (such as a  standard Markov process \cite{BG})   with state space $S$  and jointly measurable transition densities $p_{t}(x,y)$ with respect to some $\si$-finite measure $m$ on $S$.  
We assume that the potential densities  
\[
u(x,y)=\int_{0}^{\ff}p_{t}(x,y)\,dt
\]
 are   finite off the diagonal, but allow them to be infinite on the diagonal. 
  We do not require that $p_{t}(x,y)$ is  symmetric.

   We  assume furthermore that  $0<p_{t}(x,y)<\ff$ for all $0<t<\ff$ and $ x,y\in S$, and that there exists   another Borel right process $\wh X$ in duality with $X$   (Cf \cite{BG}), relative to the measure $m$, so that    its transition probabilities 
$\wh P_{t}(x,dy)=p_{t}(y,x)\,m(dy)$.  These conditions allow us to  use  material on bridge measures  in \cite{FPY}. In particular, for all $0<t<\ff$ and $x,y\in S$,   there exists a finite measure $Q_{t}^{x,y}$ on $\mathcal{F}_{t^{-}}$, of total mass $p_{t}(x,y)$, such that
\begin{equation}
Q_{t}^{x,y}\(1_{\{\ze>s\}}\,F_{s}\)=P^{x}\( F_{s}  \,p_{t-s}(X_{s},y) \),\label{10.1}
\end{equation}
for all $F_{s} \in \mathcal{F}_{s}$ with $s<t$. (We use the letter Q for measures which are not necessarily of mass $1$, and reserve the letter P for probability measures.) $Q_{t}^{x,y}$ should be thought of as a measure for paths which begin at $x$ and end at $y$ at time $t$. When normalized, this gives the bridge measure $P_{t}^{x,y}$ of \cite{FPY}.

 We use the canonical representation of $X$ in which $\Om$ is  the set of cadlag paths $\om$ in $S_{\De}=S\cup \De$ with  $\De\notin S$, and  is such that   $\om(t)=\De$ for all $t\geq \ze=\inf \{t>0\,|\om(t)=\De\}$.   Set  $X_{t}(\om)=\om(t)$.   We define a $\si$-finite measure $\mu$
on $(\Om, \mathcal{F})$ by 
\begin{equation}
\int F\,d\mu=\int_{0}^{\ff}{1\over t}\int\,Q_{t}^{x,x}\(F\circ k_{t}\)\,dm (x)\,dt\label{ls.3}
\end{equation}
for all $\mathcal{F}$ measurable functions $F$ on $\Om$.
Here $k_{t}$ is the killing operator defined by $k_{t}\om(s)=\om(s)$ if $s<t$ and $k_{t}\om(s)=\De$ if $s\geq t$, so that $k_{t}^{-1}\mathcal{F}\subset \mathcal{F}_{t^{-}}$.   As usual, if $F$ is a function, we often write $\mu (F)$ for $\int F\,d\mu$.  $\mu$ is $\si$-finite as   any set of loops  in a compact set  with lifetime bounded away from zero and infinity has finite measure.

  We call $\mu$ the loop measure of $X$ because, when $X$ has continuous paths, $\mu$ is concentrated on the set of continuous 
 loops with a distinguished starting point (since  $Q_{t}^{x,x}$ is carried by loops starting at $x$).
 Moreover, it is shift invariant. More precisely 
  let  $\rho_{u}$ denote the loop rotation defined by
\[
\rho_{u}\om (s)
	=\left\{\begin{array}{ll}
	\om (s+u \mbox{ mod } \zeta(\om)), & \mbox{ if $0\leq s<\zeta(\om)$}\\
	\De, & \mbox{otherwise.}
	\end{array} \right.
	\]
Here, for two positive numbers $a,b$ we define  $a\mbox{ mod }  b= a-mb$ as the unique positive integer $m$ such that $\,0\leq a-mb<b $ . 
 $\mu $ is invariant under  $\rho_{u}$, for any   $u$. We let $\mathcal{F}_{\rho}$ denote the $\si$-algebra of $\mathcal{F}$ measurable functions $F$ on $\Om$ which are invariant under $\rho$, that is  $F\circ \rho_{u}= F$ for all $u\geq 0$. Loop functionals of interest are mostly $\mathcal{F}_{\rho}$ - measurable.
Recall that Poisson processes of intensity $\mu$ appear naturally as produced by loop erasure in the construction of random spanning trees through the  Wilson algorithm (see chapter 8 in \cite{Le Jan1} ).   
Although the definition of  $\mu$ in (\ref{ls.3}), especially the   ${1 \over t}$, may look forbidding, $\mu$ often has a nice form when applied to specific functions in $\mathcal{F}_{\rho}$.
A particular function in $\mathcal{F}_{\rho}$ is given by
\begin{equation}
\phi(f)=\int_{0}^{\ff}f (X_{t})\,dt,\label{phi}
\end{equation}
  where $f$ is any measurable function on $S$. If $f_{j}$, $ j=1,\ldots,k\geq 2$ are non-negative functions  on $S$, then by \cite[Lemma 2.1]{LMR2}
\bea
&&
\mu\(   \prod_{j=1}^{k}\phi(f_{j})\) \label{ls.4}\\
&& = \sum_{\pi\in \mathcal{P}_{k }^{\odot}}\int u(y_{\pi(1)},y_{\pi(2)})\cdots   u(y_{\pi(k-1)},y_{\pi(k)})u(y_{\pi(k)},y_{\pi(1)})  \prod_{j=1}^{k} f_{j}(y_{j})\,dy_{j}\nn
\eea
where $\mathcal{P}_{k}^{\odot}$ denotes the set of permutations of $[1,k]$ on the circle.   (For example,    $(1,3,2)$, $(2,1,3)$ and $(3,2,1)$ are considered to be one permutation   $\pi\in \mathcal{P}^{\odot}_{3 }$.)
We note however that in general when $u$ is infinite on the diagonal
\[
\mu\(\phi(f_{j})\)=\infty. 
\]
For  $k\geq 2$, the integral $(\ref{ls.4})$ can be finite if the $f_i$ satisfy certain integrability conditions: see the beginning of section \ref{sec-multCAF}.
Consider more generally the multiple integral  
 \begin{equation}
 \sum_{\pi\in \mathcal{T}_{k}^{\odot}}  \int_{0\leq r_{1}\leq\cdots\leq r_{k}\leq t}   \,f_{\pi (1)}(X_{r_{1}})\cdots  \,f_{\pi (k)}(X_{r_{k}}). 
\end{equation}
where $\mathcal{T}_{k}^{\odot}$ denotes the set of   translations  $\pi$  of $[1,k]$ which are cyclic mod $k$, that is, for some $i$,
$\pi(j)=j+i,\mod k, $ for all  $j=1,\ldots,k$.  

Finite sums of multiple integrals such as these form an algebra (see exercise 11, p.21 in \cite{Le Jan1}) which   generates $\mathcal{F}_{\rho}$, \cite{Chang}.
 
Finally, let $ a(x) $ be a  bounded, strictly positive function on $S$.  
Define the time changed process $ Y_{t}=X_{\tau(t)}$ where $\tau(t)=\int^{t}_{0} a(Y_{s})\,ds$ is the inverse of the CAF 
 $A_{t}=\int^{t}_{0} 1/a(X_{s})\,ds$. It satisfies the duality assumption relative to the measure $a\cdot m$. It then follows as in  \cite[ section 7.3]{FR} that if $u_{X}(x,y), u_{Y}(x,y)$ denote the potential densities
  of $X,Y$ respectively with respect to $m,a\cdot m$ respectively, then
 \begin{equation}
 u_{Y}(x,y)=u_{X}(x,y)/a(y). \label{7.2}
 \end{equation}
 It follows that $\mu_{Y}$ is the image of $\mu_{X}$ by the time change.

\section{Multiple CAF's and perturbation of  loop measures}\label{sec-multCAF}

  Let $\mathcal{M}( S)$ be the set of finite signed Radon measures on $\mathcal{B}( S)$.
We  say that a norm $\|\cdot\|$ on $\MM(S)$ is a {\bf proper norm} with respect to a kernel $u$ if for all $n\geq 2$ and   $\nu_{1},\ldots, \nu_{n}$ in $\MM(S)$
  \begin{equation}
\Big | \int  \prod_{j=1}^{n}u (x_{j},x_{j+1})\prod_{i=1}^{n}\,d\nu_{j}  (x_{i} )\Big | \leq  C^{n}\prod_{j=1}^{n} \|\nu_{j}\|,\label{p.0}
 \end{equation}
  (with $x_{n+1} =x_1$) for some universal constant $C<\ff$. In Section 6 of \cite{LMR}  we present several examples  of proper norms which depend upon various hypotheses about the kernel $u$.
 
 In particular, the following norm is related to the square root of the generator of $X$, which defines the Dirichlet space in the $m$-symmetric case:
 \begin{equation}
 \|\nu \|_{w}:= \( \int\!\!\int \( \int w(x,y)w(y,z)\,d \nu(y) \)^{2}\,dm(x) \,dm(z) \)^{1/2},\label{n.lj1q}
   \end{equation}
  where   
   \begin{equation}
  w(x,y)= \int_{0}^{\ff}  {p_{s}(x,y) \over \sqrt{\pi s}}\,ds.\label{n.lj2}
   \end{equation}
   
  To see that $\|\nu\|_w$ is a proper norm we first note that 
   \begin{equation}
  u(x,z)=   \int w(x,y)w(y,z)\,dm(y) \label{n.lj3}
  \end{equation}
(It is interesting to note  that $w$ is the potential density of the process   $X_{T_{t}}$ where $T_{t}$ is the stable subordinator of index $1/2$.     In operator notation (\ref{n.lj3}) says that $W^{2}=U$ where $W$  and $U$ are   operators with kernels $w$ and $u$ respectively.)    
   
%
   Using (\ref{n.lj3})
 \bea 
    \prod_{j=1}^{n}  u  (z_{j } ,z_{j+1} )& =   &    \prod_{j=1}^{n}\int w(z_{j},\la_{j})w(\la_{j},z_{j+1}) \,dm (\la_{j})\\
    & =   & \nn  \prod_{j=1}^{n}\int w(z_{j},\la_{j})w(\la_{j-1},z_{j }) \,dm (\la_{j})
   \eea
in which    $z_{n+1}=z_{1}$ and $\la_{0}=\la_{n}$. 
  It follows from this that  
  \bea
\lefteqn{
\Big | \int  \prod_{j=1}^{n}u(z_{j } ,z_{j+1} )\prod_{j=1}^{n}\,d\nu_{j}  (z_{j} )\Big |\label{4.7a}}\\   &= &  \Bigg| \int \prod_{j=1}^{n} \( \int w(z_{j},\la_{j})w(\la_{j-1},z_{j}) \,d\nu_{j} (z_{j})
 \)\prod_{j=1}^{n}\,dm (\la_{j})\Bigg
 |\nn \\
  &\le &\prod_{j=1}^{n}\( \int\!\!\int \( \int w(z_{j},s)w(t,z_{j}) \,d\nu_{j} (z_{j})\)^{2}\,dm (s)\,dm (t)\)^{1/2}\leq \prod_{j=1}^{n} \|\nu_{j}\|_{w},\nn
  \eea 
where, for the first inequality, we use repeatedly   the Cauchy-Schwarz inequality.

Lastly, set 
\be
   M_{\nu}(x,z)=\int w(x,y)w(y,z)\,d \nu(y).\label{6.16}
\ee
  Since $\|\nu\|_w$ is the $L^{2}$ norm of $M_{\nu}$, and  $M_{\nu+\nu'}=M_{\nu}+M_{\nu'}$, we see that  $\|\nu\|_w$ is a  norm. (This can also be viewed as the Hilbert-Schmidt norm of the operator  defined by the kernel $M_\nu$).
  \qed

   We denote by  $\mathcal{R}^{+}$  the set of positive bounded  Revuz  measures $\nu$ on $S$   that are associated with $X$. This is explained in detail in  Section 2.1 of \cite{LMR}. We use $L_{t}^{\nu}$ to denote the CAF with Revuz measure $\nu$.

Let  $\|\cdot\|$ be a  proper norm on   $\MM(S)$ with respect to the kernel  $u $. Set  
\begin{equation}
\MM^{+}_{\|\cdot \|}=\{\mbox{positive }\nu \in \MM (S)\,|\, \|\nu \|<\ff\},\label{con.1i}
\end{equation} 
and
 \begin{equation}
\mathcal{R}^{+}_{\|\cdot\|}=\mathcal{R^{+}}\cap \MM^{+}_{\|\cdot\|}.\label{rev1i}
 \end{equation}
Let  $\MM_{\|\cdot \|}$ and $\mathcal{R}_{\|\cdot\|}$  denote the set of measures  of the form $\nu=\nu_{1}-\nu_{2}$ with $ \nu_{1},\nu_{2}\in  \MM^{+}_{\|\cdot \|}$ or $ \mathcal{R}^{+}_{\|\cdot\|}$ respectively.  
We   often omit saying that both $\RR_{\|\cdot\|}$ and $\|\cdot\|$ depend on the kernel $u$.

Let $\|\cdot\|$ be a proper  norm for $u$. For $\nu_{j}\in  \mathcal{R}_{\|\cdot\|}  $, $j=1,\ldots,k$,  let 
\begin{equation}
M^{\nu_{1},\ldots,\nu_{k}}_{t}= \sum_{\pi\in  \mathcal{T}_{k}^{\odot}}  \int_{0\leq r_{1}\leq\cdots\leq r_{k}\leq t}   \,dL^{\nu_{\pi (1)}}_{r_{1}}\cdots  \,dL^{\nu_{\pi (k)}}_{r_{k}}. \label{2.0}
\end{equation}
 We refer to $M^{\nu_{1},\ldots,\nu_{k}}_{t}$ as a multiple CAF.  
 Let
\begin{equation}
Q^{x,y} \(F\)=\int_{0}^{\ff}  \,Q_{t}^{x,y}\(F\circ k_{t}\)\,dt.\label{Q.3}
\end{equation}
 We have the following analogue of \cite[Proposition 5]{Le Jan1} and \cite[Lemma  2.1]{LMR}.

\begin{lemma}\label{lem-mmomm}  For any measures  $\nu_{j}\in \mathcal{R}^{+}_{\|\cdot\|}  $, $j=1,\ldots, k\geq 2$,  
\begin{equation}
\mu(  M^{\nu_{1},\ldots,\nu_{k} }_{\ff})= \int u(y_{1},y_{2})\cdots   u(y_{k-1},y_{k})u(y_{k},y_{1})  \prod_{j=1}^{k} \,d\nu_{j}(y_{j}),\label{2.1}
\end{equation}
\bea
&&
\hspace{-.3 in} Q^{x,y} (  M^{\nu_{1},\ldots,\nu_{k} }_{\ff}) \label{2.1a}\\
&&\hspace{-.3 in} =\sum_{\pi\in  \mathcal{T}_{k}^{\odot}}\int u(x,y_{1})u(y_{1},y_{2})\cdots   u(y_{k-1},y_{k})u(y_{k},y)  \prod_{j=1}^{k} \,d\nu_{\pi(j)}(y_{j}),\nn
\eea
and if $\nu \in  \mathcal{R}^{+}_{\|\cdot\|}  $
\begin{eqnarray}
&&\mu(  M^{\nu_{1},\ldots,\nu_{k} }_{\ff}\,L^{\nu }_{\ff})
\label{2.1b}\\
&&   =\sum_{i=1}^{k}\int\(\prod_{j=1}^{i-1}u (y_{j},y_{j+1})\)\nn\\
&& \hspace{1 in} u (y_{i},x) u(x,y_{i+1})  \(  \prod_{j=i+1}^{k}u (y_{j},y_{j+1})\)  \prod_{j=1}^{k} \,d\nu_{j}(y_{j})\,d\nu(x)  \nonumber\nn\\
&&=\int Q^{x,x} (  M^{\nu_{1},\ldots,\nu_{k} }_{\ff})\,d\nu(x),  \nonumber
\end{eqnarray}
with $y_{k+1}=y_{1}$.
\end{lemma}

The proof of (\ref{2.1}) follows that of \cite[Lemma  2.1]{LMR}, noticing that the crucial step \cite[(2.23)-(2.28)]{LMR} used the fact that the set of permutations of $[1,k]$ is invariant under translation mod $k$. Since 
$\mathcal{T}_{k}^{\odot}$ is invariant under translation mod $k$, the same proof will work here. The proof of (\ref{2.1a}), which is much easier, follows that of \cite[Lemma  4.2]{LMR}.
The first equality in (\ref{2.1b}) follows from (\ref{2.1}) and the fact that
\begin{equation}
M^{\nu_{1},\ldots,\nu_{k} }_{\ff}\,L^{\nu }_{\ff}=\sum_{i=1}^{k} M^{\nu_{1},\ldots,\nu_{i-1},\nu,\nu_{i},\ldots,\nu_{k} }_{\ff},\label{2.1c}
\end{equation}
which is easily verified using $L^{\nu }_{\ff}=\int_{0}^{\ff}\,dL^{\nu }_{t}$. The second equality in (\ref{2.1b}) then follows by comparing with (\ref{2.1a}).
\qed

Let now $X_{(\ep)}$, $\ep\geq 0$, $X_{(0)}=X$, be a family of Markov processes  with potential densities $u_{(\ep)}(x,y)$, and let       $\mu_{(\ep)}$ denote the loop measure for $X_{(\ep)}$. Assume that we can use the same proper norm $\|\cdot\|$ for all $u_{(\ep)}$. Let
\begin{equation}
u_{(0)}'(x,y)={du_{(\ep)}(x,y) \over d\ep}\,\,|_{\ep=0},\label{2.3}
\end{equation} 
and assume that $\|\cdot\|$ is also a proper norm for $u_{(0)}'$.
Then using the last Lemma we have formally, that is, assuming we can justify interchanging  
derivative and integral in the second equality,  
 \begin{eqnarray}
&&{d \over d\ep}\,\,\,\mu_{(\ep)}( M^{\nu_{1},\ldots,\nu_{k} }_{\ff})\Big |_{\ep=0}
\label{2.4}\\
&& ={d \over d\ep}\,\,\,\int u_{(\ep)}(y_{1},y_{2})\cdots  u_{(\ep)}(y_{k-1},y_{k})u_{(\ep)}(y_{k},y_{1})  \prod_{j=1}^{k} \,d\nu_{j}(y_{j})\Big |_{\ep=0}\nn  \nonumber\\
&& = \sum_{i=1}^{k}\int u(y_{1},y_{2})\cdots   u_{(0)}'(y_{i},y_{i+1})\cdots u(y_{k},y_{1})  \prod_{j=1}^{k} \,d\nu_{j}(y_{j})  \nonumber\\
&& =\sum_{\pi\in \mathcal{T}_{k}^{\odot}}
 \int  u(y_{1},y_{2})\cdots   u(y_{k-1},y_{k })u_{(0)}'(y_{k},y_{1}) \prod_{j=1}^{k} \,d\nu_{\pi(j)}(y_{j})  \nonumber
\end{eqnarray}
where we have set $y_{k+1}=y_{1}$.

Assume now that for some distribution $F$ on $S\times S$ we have 
\begin{equation}
u_{(0)}'(y_{k},y_{1})=\int_{S\times S}u(y_{k },x)F(x,y) \,u(y,y_{1})\,dm(x)\,dm(y).\label{2.5}
\end{equation}
Let $\AA_{\|\cdot\|}$ denote the space spanned by the multiple CAF's with $\nu_{j}\in \RR_{\|\cdot\|}$. Note it is an algebra. 
Then comparing (\ref{2.4}) and (\ref{2.1a}) we would obtain
\begin{equation}
{d\mu_{(\ep)}(A) \over d\ep} \,\,|_{\ep=0}=\int_{S\times S}F(x,y) Q^{ y,x}(A)\,dm(x)\,dm(y),\label{2.6}
\end{equation}
for all $A\in\AA_{\|\cdot\|}$.

In the following sections we present specific examples where  this heuristic approach is made rigorous.

\section{Perturbation of  L\'evy processes}\label{sec-plev}

Let $X$ be a transient L\'evy process in $R^{d}$ with characteristic exponent $\psi$ so that, as distributions
\begin{equation}
u(x,y)=\int {e^{i\la (y-x)} \over \psi(\la)}\,d\la.\label{4.1}
\end{equation}
In \cite{LMR} we showed that $\| \cdot \|_{\psi, 2}$ is a proper norm for $u$ where
\begin{equation}
\|  \nu \|_{\psi, 2}^{2}=\int \({1 \over |\psi |}\ast {1 \over |\psi |}\,(\la)\)\,\,|\wh \nu (\la)|^{2}\,d\la.\label{3.7}
\end{equation}

Let $\ka$  be a L\'evy  characteristic exponent, so that the same is true for $\psi+\ep\ka$, and let $X_{(\ep)}$ be the  L\'evy process with characteristic exponent $\psi+\ep\ka$. We let $u_{(\ep)}(x,y)$ denote the potential of $X_{(\ep)}$ so that, as distributions
\begin{equation}
u_{(\ep)}(x,y)=\int {e^{i\la (y-x)} \over \psi(\la)+\ep\ka (\la)}\,d\la.\label{4.2}
\end{equation}
If we assume that 
\begin{equation}
| \ka (\la)|\leq C |\psi(\la)|\label{3.8}
\end{equation}
for some $C<\ff$, then for $\ep >0$ sufficiently small
\begin{equation}
\|  \nu \|_{\psi+\ep\ka, 2}\leq C'  \|  \nu \|_{\psi, 2},\label{3.9}
\end{equation}
for some $C'<\ff$.
Thus $\| \cdot \|_{\psi, 2}$ is a proper norm for $u_{(\ep)}$.

Let $F$ be the distribution   given by
\begin{equation}
F(x,y)=\int  e^{i\la (y-x)}   \ka(\la)\,d\la,\label{4.3}
\end{equation}
and let $\wh Q^{\la_{1},\la_{2}}(A)$  denote the Fourier transform of $Q^{x,y}(A)$ in $x,y$.
\bt\label{prop-levy} If (\ref{3.8}) holds, then
 \bea
{d\mu_{(\ep)}(A) \over d\ep}\,\,\,|_{\ep=0}&=&-\int_{R^{d}\times R^{d}}  Q^{y,x}(A)F(x,y) \,dm(x)\,dm(y)\label{3.6rf}\\
&=&  -\int  \wh Q^{\la ,-\la }(A)\,\,   \ka(\la)\,d\la.  \nn
\eea
for all $A\in\AA_{\|\cdot\|_{\psi, 2}}$.
\et

{\bf  Proof of Theorem \ref{prop-levy}: }It suffices to show that for any  $ \nu_{1},\ldots, \nu_{k}\in \mathcal{R}_{\|\cdot\|_{\psi, 2}}^{+} $
\begin{eqnarray}
&&{d \over d\ep}\int \prod_{j=1}^{k}u_{(\ep)}(y_{j},y_{j+1})  \prod_{j=1}^{k} \,d\nu_{j}(y_{j})\,\,\,\Big|_{\ep=0}
\label{le.4}\\
&&   =-\sum_{i=1}^{k}\int\(\prod_{j=1}^{i-1}u (y_{j},y_{j+1})\)\nn\\
&& \hspace{1.3 in}\(\int_{R^{d}\times R^{d}}    u(y_{i},x)\,F(x,y)\, u(y,y_{i+1}) \,dm(x)\,dm(y)\)\nn\\
&& \hspace{2.3 in}    \(  \prod_{j=i+1}^{k}u (y_{j},y_{j+1})\)  \prod_{j=1}^{k} \,d\nu_{j}(y_{j}),  \nonumber
\end{eqnarray}
with $y_{k+1}=y_{1}$.

Using (\ref{4.2}) we see that  
 \bea
I(\ep)&=:& \int  \prod_{j=1}^{k}u_{(\ep)}(y_{j},y_{j+1}) \prod_{j=1}^{k}\,d\nu_{j}  (y_{j} ) \label{m8.3q}\\
&=&     \int\!\! \int  \prod_{j=1}^{k} e^{i(y_{j+1}-y_{j}) \cd \la_{j}}{1 \over \psi(\la_{j})+\ep\ka (\la_{j})}\,d\la_{j} d\nu_{j}(y_{j})  \nn \\
   &=&    \int \( \prod_{j=1}^{k}  \int e^{ -i(\la_{j }-\la_{j-1}) \cd y_{j}}
\,d\nu_{j}(y_{j}) \) \prod_{j=1}^{k} {1 \over \psi(\la_{j})+\ep\ka (\la_{j})}\,d\la_{j}   \nn
  \eea
   where $\la_{0}=\la_{n}$.  
   We take the Fourier transforms of the $\{\nu_{j}\}$ to see that 
   \bea
&& \int  \prod_{j=1}^{k}u_{(\ep)}(y_{j},y_{j+1}) \prod_{j=1}^{k}\,d\nu_{j}  (y_{j} ) \label{le.8}\\
&&  =   \int  \prod_{j=1}^{k}  \hat\nu_{j}(\la_{j }-\la_{j-1}) {1 \over \psi(\la_{j})+\ep\ka (\la_{j})} \,d\la_{j}.\nn
  \eea
  We have
  \begin{equation}
{1 \over \psi(\la_{j})+\ep\ka (\la_{j})}={1 \over \psi(\la_{j})}-\ep {\ka (\la_{j}) \over \psi^{2}(\la_{j})}+\ep^{2} {\ka^{2} (\la_{j}) \over \psi^{2}(\la_{j})(\psi(\la_{j})+\ep\ka  (\la_{j}))}.  \label{le.9}
  \end{equation}
  Substituting this into (\ref{le.9}) we can write the result as 
  \begin{equation}
I(\ep)=I(0)-\ep J +K(\ep),  \label{le.10}
  \end{equation}
   where
   \begin{equation}
  J= \sum_{i=1}^{k} \int \(\prod_{\stackrel{j=1}{j\neq i}}^{k}{1 \over \psi(\la_{j}) }\) {\ka (\la_{i}) \over \psi^{2}(\la_{i})}  \prod_{j=1}^{k}  \hat\nu_{j}(\la_{j }-\la_{j-1})  \,d\la_{j},\label{le.11}
   \end{equation}
   and $K(\ep)$ is the sum of all remaining terms.
   
   We now show that $J$ is precisely the right hand side of (\ref{le.4}), and that 
   \begin{equation}
|K(\ep)|=O(\ep^{2}),\label{le.12}
\end{equation}
which will complete the proof of our Proposition. 

For the first point, using the relation between Fourier transforms and convolutions we have
\bea
&&
 \int_{R^{d}\times R^{d}}       u(y_{i},x)\,F(x,y)\, u(y,y_{i+1}) \,dm(x)\,dm(y)\label{le.13}\\
 &&\hspace{2 in}
= \int e^{i\la_{i} (y_{i+1}-y_{i})}{ \ka (\la_{i} )\over \psi^{2}(\la_{i} ) }\,d\la_{i}. \nn
\eea 
Using this in the right hand side of (\ref{le.4}) and proceeding as in (\ref{m8.3q})-(\ref{le.8}) we indeed obtain $J$. As for (\ref{le.12}), $K(\ep)$ is the sum many terms each of which has a factor of $\ep^{m}$ for some $m\geq 2$. We need only show that the corresponding integrals are bounded uniformly in $\ep$. For example, consider the term which arises when using the last term in (\ref{le.9}) for all $j$. This term is $\ep^{2k}$ times
   \be  \wt K(\ep)=: \int  \prod_{j=1}^{k}  \hat\nu_{j}(\la_{j }-\la_{j-1}){\ka^{2} (\la_{j}) \over \psi^{2}(\la_{j})(\psi(\la_{j})+\ep\ka  (\la_{j}))} \,d\la_{j}. \label{le.14}
  \ee 
By our assumption (\ref{3.8}), for sufficiently small $\ep$
\bea
| \wt K(\ep)| &\leq &C'\int \prod_{j=1}^{k} | \hat\nu_{j}(\la_{j }-\la_{j-1})|{1\over |\psi(\la_{j})|} \,d\la_{j} \label{le.15}\\
&\leq &C''\prod_{j=1}^{k}\|  \nu_{j} \|_{\psi, 2},\nn
\eea
by repeated use of the Cauchy-Schwarz inequality, as in our proof in \cite{LMR}   that $\| \cdot \|_{\psi, 2}$ is a proper norm for $u$.
\qed

\section{Perturbation by multiplicative functionals}\label{sec-pmult}

Let $m_{t}$ be a continuous decreasing multiplicative functional of $X$, with $m_{t}\leq 1$ for all $t$ and $m_{\ze}=0$. By \cite[Theorem 61.5]{S},
there is a right process $\wt{X}_{t}$ with transition semigroup  
\begin{equation}
\wt{P}_{t}f(x)=P^{x}(f(X_{t})m_{t})=\int Q_{t}^{x,y}(m_{t})f(y)\,dm(y),\label{m.1}
\end{equation}
where $Q_{t}^{x,y}$ are the bridge measures (\ref{10.1}) for our original process $X$ and the second equality is \cite[(2.8)]{FPY}. Thus $\wt{X} $ has transition densities
\begin{equation}
\wt p_{t}(x,y)=: Q_{t}^{x,y}(m_{t}),\label{m.2}
\end{equation}
and using \cite[(2.8)]{FPY} once more we can verify that these satisfy the Chapman-Kolmogorov  equations. It follows from the construction in \cite[Theorem 61.5]{S} that if $X$ has cadlag paths  so will $\wt{X}$. Let    $r_{t}(\om )(s)=\om (t-s)^{-}$, $0\leq s\leq t$ be the time reversal mapping. 
If we set $\wh m_{t}=m_{t}\circ r_{t}$, then it is easy to check that $\wh m_{t}$ is a multiplicative functional as above, see \cite[p. 359]{CW}. If $\wh X$ is the dual process for $X$ as described in Section 2, with bridge measures $\wh Q_{t}^{x,y}$ then as above there exists a process $Y$ with transition densities 
$\wh Q_{t}^{x,y}(\wh m_{t})$. By \cite[Corollary 1]{FPY}
\begin{equation}
\wh Q_{t}^{x,y}(\wh m_{t})=\wh Q_{t}^{x,y}( m_{t}\circ r_{t})=Q_{t}^{y,x}( m_{t}),\label{m.2q}
\end{equation}
which shows that $Y$ is dual to $\wt{X}$.

We now show that  if $\wt Q_{t}^{x,y}$ are the bridge measures for $\wt{X} $, then
\begin{equation}
\wt Q_{t}^{x,y}(F)=Q_{t}^{x,y}(F\,m_{t}), \mbox{ for }F\in \mathcal{F}_{s}, \,s<t.\label{m.3}
\end{equation}
To see this, using the fact that $m_{t}$ is continuous and  decreasing, we have for  $F\in \mathcal{F}_{s}, \,s<t$
\begin{eqnarray}
\wt Q_{t}^{x,y}(F)&=&\wt{P}^{x}\(F Q_{t-s}^{X_{s},y}(m_{t-s})\)
\label{m.4}\\
&=&\lim_{t^{\ast}\uparrow t} \wt{P}^{x}\(F Q_{t-s}^{X_{s},y}(m_{t^{\ast}-s})\)  \nonumber\\
&=&\lim_{t^{\ast}\uparrow t} \wt{P}^{x}\(F P^{X_{s} }\(m_{t^{\ast}-s}\,p_{t-t^{\ast}}(X_{t^{\ast}-s},y)      \)\)  \nonumber\\
&=&\lim_{t^{\ast}\uparrow t} P^{x}\(F m_{s} P^{X_{s} }\(m_{t^{\ast}-s}\,p_{t-t^{\ast}}(X_{t^{\ast}-s},y)      \)\)  \nonumber\\
&=&\lim_{t^{\ast}\uparrow t} P^{x}\(F m_{s}\,m_{t^{\ast}-s}\circ\th_{s}  \,p_{t-t^{\ast}}(X_{t^{\ast}},y)      \)  \nonumber\\
&=&\lim_{t^{\ast}\uparrow t} P^{x}\(F m_{t^{\ast}}  \,p_{t-t^{\ast}}(X_{t^{\ast}},y)      \)  \nonumber\\
&=&\lim_{t^{\ast}\uparrow t} Q_{t}^{x,y}(F\,m_{t^{\ast}}) =Q_{t}^{x,y}(F\,m_{t}), \nonumber
\end{eqnarray}
which proves (\ref{m.3}). 

If $A_{t}$ is a CAF,  then $m_{t}=e^{-A_{t}}$ is a continuous decreasing multiplicative functional of $X$.   Let $X_{(\ep)}$ denote the Markov process $\wt{X}$ with $m_{t}=e^{-\ep L^{\nu}_{t}}$.
If $\mu_{(\ep)}$ denotes the loop measure for $X_{(\ep)}$, and $\mu$ the loop measure for $X$, it now follows from (\ref{ls.3}) and (\ref{m.3}) that
\begin{equation}
\mu_{(\ep)}\(A\)=\mu \(Ae^{-\ep L^{\nu}_{\ff}}\).\label{m.5}
\end{equation}

\bt\label{prop-time}
If $\nu\in  \mathcal{R}_{\|\cdot\|}^{+} $, for some proper norm $\|\cdot\|$, then
\begin{equation}
{d\mu_{(\ep)}(A) \over d\ep}\,\,\,|_{\ep=0}=- \mu (L_{\ff}^{\nu}A)=-\int_{S}  Q^{x,x}(A)\,d\nu(x) \label{3.6}
\end{equation}
for all $A\in\AA_{\|\cdot\|}$.
\et

{\bf  Proof: } It suffices to prove this for $A$ of the form $M^{\nu_{1},\ldots,\nu_{k} }_{\ff}$ with $\nu_{j}\in  \mathcal{R}_{\|\cdot\|}^{+}$, $ 1\leq j\leq k $. Since $0\leq e^{-x}-1+x\leq x^{2}/2$ for $x\geq 0$, it follows from (\ref{m.5}) that 
\begin{equation}
|\mu_{(\ep)}(M^{\nu_{1},\ldots,\nu_{k} }_{\ff})-\mu (M^{\nu_{1},\ldots,\nu_{k} }_{\ff})-\ep \mu (L^{\nu}_{\ff}M^{\nu_{1},\ldots,\nu_{k} }_{\ff})|\leq \ep^{2} \mu \(\(L^{\nu}_{\ff}\)^{2}M^{\nu_{1},\ldots,\nu_{k} }_{\ff}\).\label{m.7}
\end{equation}
$\mu \(\(L^{\nu}_{\ff}\)^{2}A\)$ is bounded by our assumption about the proper norm $\|\cdot\|$, so the first equality in (\ref{3.6}) follows. The second equality  is (\ref{2.1b}).\qed

It follows from (\ref{m.2}) that $X_{(\ep)}$ has potential densities
\begin{equation}
 u_{(\ep)}(x,y)= \int_{0}^{\ff}Q_{t}^{x,y}(e^{-\ep L^{\nu}_{t}})\,dt\leq u(x,y).\label{pd.1}
\end{equation}
Then, using \cite[Lemma 1]{FPY}
\begin{eqnarray}
 u'_{(0)}(x,y)&=& -\int_{0}^{\ff}Q_{t}^{x,y}( L^{\nu}_{t})\,dt
\label{pd.2}\\
&=& -\int_{0}^{\ff}\int_{0}^{t} \int p_{s}(x,z)p_{t-s}(z,y) \,d\nu(z)\,ds\,dt  \nonumber\\
&=& -\int u(x,z)u(z,y)\,d\nu(z). \nonumber
\end{eqnarray}
In view of this, Theorem \ref{prop-time} is another example of our heuristic formula (\ref{2.6}), where the distribution $F$ on $S\times S$ of (\ref{2.5}) is $\de (x-y)\,d\nu(x)$. 

For use in the next section we now recast Theorem \ref{prop-time}. 
For $\nu_{j}\in  \mathcal{R}_{\|\cdot\|}^{+}$, $ 1\leq j\leq k $ let
\begin{equation}
I(\ep)=\int \prod_{j=1}^{k}u_{(\ep)}(y_{j},y_{j+1})  \prod_{j=1}^{k} \,d\nu_{j}(y_{j}).\label{rg.1}
\end{equation}
Since $\int \prod_{j=1}^{k}u_{(\ep)}(y_{j},y_{j+1})  \prod_{j=1}^{k} \,d\nu_{j}(y_{j})=\mu_{(\ep)}(M^{\nu_{1},\ldots,\nu_{k} }_{\ff})$ by (\ref{2.1}), it follows from Theorem \ref{prop-time} and (\ref{2.1b})
  that 
\begin{eqnarray}
&&\lim_{\ep\to 0}{I(\ep)-I(0) \over \ep}
\label{df.7}\\
&&=  -\sum_{i=1}^{k}\int\(\prod_{j=1}^{i-1}u (y_{j},y_{j+1})\) u(y_{i},x)u(x,y_{i+1}) \nn\\
&& \hspace{2.3 in}    \(  \prod_{j=i+1}^{k}u (y_{j},y_{j+1})\)  \prod_{j=1}^{k} \,d\nu_{j}(y_{j})\,d\nu (x).  \nonumber
\end{eqnarray}

In the next section we will also use the following observation. If $m_{t}=e^{-\ep \int_{0}^{t}c(s)\,ds}$ then $\wh m_{t}=m_{t}\circ r_{t}=m_{t}$. Hence if $\wh X$ is the dual process for $X$ as described in Section 2, it follows from (\ref{m.2q}) that $Y$, the dual  process to  $X_{(\ep)}$, is the process  $(\wh X)_{(\ep)}$ obtained by perturbing $\wh X$ by the same multiplicative functional $m_{t}$. In particular
\begin{equation}
\wh u_{(\ep)}(x,y)=u_{(\ep)}(y,x).\label{swi}
\end{equation}

\section{Perturbation by addition of jumps}\label{sec-jump}

 Let $j(x,y)$ be a nonnegative $m\otimes m$-integrable function on $S\times S$.  We will assume that 
 \begin{equation}
  c(x)=:\int j(x,y)m(dy)=\int j(y,x)m(dy).\label{jl.1}
 \end{equation}
$  c(x)$ is integrable and we will assume moreover that it is bounded and strictly positive. Then
 \begin{equation}
 G(x,dy)=:\frac{1}{c(x)}j(x,y)m(dy)\label{jl.1g}
 \end{equation}
 is a probability kernel on $S\times \BB (S)$.    $c(x)$ will govern the rate at which we will add jumps  to the process, which may depend on the position $x$ of the process, and $G(x,dy)$ will describe the distribution of the jumps from position $x$. 

 In more detail, define the CAF
\begin{equation}
A_{t}=\int_{0}^{t}  c(X_{s})\,ds,\label{5.1}
\end{equation}
and let $\tau_{t}$ be the right continuous inverse of $A_{t}$. Let  $\la$ be an independent mean $1$ exponential. We define a new process $Y_{t}$ to be equal to $X_{t}$ for $t<\tau_{\la}$, and then re-birthed at a random point independent of $\lambda$, distributed according to  $G(X_{\tau_{\la}},dy) $, with this process being iterated. We use $U_{c,G}$ to denote the potential operator of $Y$.

Let
\begin{equation}
\|\nu\|_{u^2,\ff}:= |\nu|(S)\vee\,\sup_{x}\int \( u^{2}(x,y) + u^{2}(y,x) \)\,d|\nu|  (y),\label{p.2}
\end{equation}
 where $|\nu|$ is the total variation of the measure $\nu$. This is a proper norm for $u$, see \cite[(3.25)]{LMR}.
 
  We use $\mu_{(\ep)}$ to denote the loop measure associated to $Y$, where we have replaced $c$ by $\ep c$.

 \begin{theorem}\label{theo-jump}
  Assume that 
 \begin{equation}
\sup_{z}\int   u ( z , y) dm(y)   <\ff,\hspace{.2 in}\sup_{z}\int   u ( y , z) dm(y)   <\ff.\label{gc.4}
 \end{equation}
 Then $\mu_{(\ep)}$ is well defined for $\epsilon$ small enough and
  \begin{equation}
{d\mu_{(\ep)}(A) \over d\ep} \,\,|_{\ep=0}=\int_{S\times S} \(Q^{ y,x}(A)-Q^{x,x}(A)\) c(x) G(x,\,dy)\,dm(x),\label{5.10}
\end{equation}
for all $A\in\AA_{\|\cdot\|_{u^2,\ff}}$.
  \end{theorem}
  
  Before proving this theorem, we first show that  $U_{\ep c,G}$ has densities $u_{\ep c,G}(x,y)$ for $\ep$ sufficiently small.
    Note first that   by (\ref{jl.1})-(\ref{jl.1g}) 
    \begin{equation}
  \int  c(z)G(z,dy)\, dm(z)=  \int  j(z,y)\, dm(y) \, dm(z)=c(y)\, dm(y),  \label{jl.1gr}
    \end{equation}
and    since we assumed that  $c$ is bounded   it follows from   (\ref{gc.4}) that
  \begin{equation}
\sup_{x}\int  c(z)G(z,dy)u  ( y , x)\, dm(z)=\sup_{x}\int  c(y)u  ( y , x)\,  dm(y) <\ff.\label{gc.4f}
 \end{equation}
 
Let  $\la_{1}, \la_{2},\ldots $ be a sequence of independent mean $1$ exponentials, and set $T_{j}=\sum_{i=1}^{j}\la_{j}$.  
Using the fact that
 $\tau_{t+u}=\tau_{t}+\tau_{u}\circ \th_{\tau_{t}}$
 we see that
 \begin{eqnarray}
W_{c,G,n}f(x) &=:& P^{x}\(\int_{\tau(T_{n})}^{\tau(T_{n+1})}f\(Y_{t}\)\,dt\) 
 \label{gc.10}\\
 &=&P^{x}\(\(\int_{0}^{\tau(\la_{n+1})}f\(Y_{t}\)\,dt\)\circ \th_{\tau (T_{n})} \)   \nonumber\\
&=&P^{x}\(P^{Y_{\tau(T_{n})}}\(\int_{0}^{\tau(\la_{n+1})}f\(X_{t}\)\,dt\)  \)   \nonumber\\
&=&P^{x}\(\int G(Y^{-}_{\tau(T_{n})},\,dz )P^{z }\(\int_{0}^{\tau(\la_{n+1})}f\(X_{t}\)\,dt\)  \).   \nonumber
 \end{eqnarray}
 We have 
\begin{eqnarray}
&&P_{\la}^{x}\( \int_{0}^{\tau_{\la}}f\(X_{t}\) \,dt\)
\label{5.3}\\
&& =P_{\la}^{x}\( \int_{0}^{\ff}1_{\{t<\tau_{\la}\}}f\(X_{t}\) \,dt\)  \nonumber\\
&& =P_{\la}^{x}\( \int_{0}^{\ff}1_{\{A_{t}< \la\}}f\(X_{t}\) \,dt\)  \nonumber\\
&& =P^{x}\( \int_{0}^{\ff}e^{-A_{t}}f\(X_{t}\) \,dt\).  \nonumber
\end{eqnarray}
Hence setting  
\begin{equation}
V_{c}f(x)=P^{x}\( \int_{0}^{\ff}e^{-A_{t}}f\(X_{t}\) \,dt\),\label{5.6}
\end{equation}
and writing $Gh(x)=\int_{S}G(z,dy)h(y)$ for any nonnegative or bounded function $h$, we have shown that
  \begin{eqnarray}
W_{c,G,n}f(x) &=& P^{x}\(\int G(Y^{-}_{\tau(T_{n})},\,dz )V_{c}f (z )   \) 
 \label{gc.11}\\ 
 &=&P^{x}\(  GV_{c}f (Y^{-}_{\tau(T_{n})})   \).   \nonumber
 \end{eqnarray}
 
Using once again the fact that
 $\tau_{t+u}=\tau_{t}+\tau_{u}\circ \th_{\tau_{t}}$ and the Markov property, we see that for any $h$ 
 
 \begin{eqnarray}
 P^{x}\(  h(Y^{-}_{\tau(T_{n})})   \)&=&P^{x}\(  \(h(Y^{-}_{\tau(\la_{n})}) \)\circ \th_{\tau (T_{n-1})}  \)
 \label{gc.12}\\
&=& P^{x}\( P^{Y_{\tau(T_{n-1})}} \(h(Y^{-}_{\tau(\la_{n})}) \)  \)  \nonumber\\
&=& P^{x}\( \int G(Y^{-}_{\tau(T_{n-1})},\,dz ) P^{z } \(h(X^{-}_{\tau(\la_{n})}) \)  \).  \nonumber
 \end{eqnarray}
  Using the change of variables formula, \cite[Chapter 6, (55.1)]{DM2}, and the fact that $X^{-}_{t}=X_{t}$ for a.e. t, we see that 
\bea
P_{\la}^{z}\(h(X^{-}_{\tau_{\la}}) \)=P_{\la}^{z}\(h(X_{\tau_{\la}}) \)&=&E^{x}\(\int_{0}^{\ff}e^{-t} h\(X_{\tau_{t}}\)\,dt\)\label{6.82adn}\\
&=&
P^{z}\(\int_{0}^{A_{\ze}}e^{-t} h\(X_{\tau_{t}}\)\,dt\)\nn\\
&=&P^{z}\(\int_{0}^{\ff}e^{-A_{s}} h\(X_{s}\)\,dA_{s}\)\nn\\
&=&P^{z}\(\int_{0}^{\ff}e^{-A_{s}} h\(X_{s}\)c(X_{s})\,ds\)\nn\\
&=&   V_{c}(ch)(z).  \nn
\eea
Thus we can write (\ref{gc.12}) as
 \begin{equation}
 P^{x}\(  h(Y^{-}_{\tau(T_{n})})   \)=P^{x}\(  GV_{c}\,c\,h (Y^{-}_{\tau(T_{n-1})})   \). \label{gc.13}
 \end{equation}
Iterating this we obtain
\begin{eqnarray}
 P^{x}\(  h(Y^{-}_{\tau(T_{n})})   \)&=&P^{x}\(  GV_{c}\,c\,h(Y^{-}_{\tau(T_{n-1})})   \)
\label{gc.13d}\\
&=&P^{x}\(  (GV_{c}\,c)^{2}\,h (Y^{-}_{\tau(T_{n-2})})   \)  \nonumber\\
&=&\cdots  \nonumber\\
&=&P^{x}\(  (GV_{c}\,c)^{n-1}\,h (Y^{-}_{\tau(\la)})   \)  \nonumber\\
&=&V_{c}\,c (GV_{c}\,c)^{n-1}\,h (x),      \nonumber
\end{eqnarray}
where the last step used (\ref{6.82adn}).
 Applying this to (\ref{gc.11})  we have that
 \begin{equation}
W_{c,G,n}f(x)=V_{c}\,c (GV_{c}\,c)^{n-1}\,GV_{c}f (x) =V_{c}\(c GV_{c}\)^{n} f(x).\label{gc.14}
 \end{equation}
 It follows from  (\ref{pd.1}) that $V_{ c}$ has a density which we write as $v_{ c}(x,y)$, and therefore $W_{c,G,n}$ has the density
    \bea
 &&
\hspace{-.3 in}w_{c,G, n} (x,y)  =\int v_{  c}(x,z_{1})G(z_{1},\,dz_{2})v_{  c}( z_{2}, z_{3})\cdots \label{gc.20}\\
&&\hspace{-.1 in}\cdots G(z_{2n-3},\,dz_{2n-2})v_{  c}( z_{2n-2}, z_{2n-1})G(z_{2n-1},\,dz_{2n})v_{  c}( z_{2n}, y)\nn\\
&&\hspace{3 in}\prod_{j=1}^{n}c(z_{2j-1})\,dm(z_{2j-1}).\nn
 \eea

By (\ref{gc.4}) it follows from  (\ref{pd.1})  that  
  \begin{equation}
\sup_{z}\int   v_{ c} ( z , y) dm(y) \leq M,\label{gc.4b}
 \end{equation}
 and thus by (\ref{gc.20}) that
   \begin{equation}
\sup_{z}\int   w_{c,G, n}  ( z , y) dm(y) \leq C^{n}M^{n+1}.\label{gc.4cma}
 \end{equation}

 Replacing $c$ by $\ep c$ we have shown that for $\ep$ sufficiently small
 \begin{equation}
U_{\ep c, G}f(x)= \sum_{n=0}^{\ff}\ep^{n }\,\,V_{\ep c}\(c GV_{\ep c}\)^{n} f(x)\label{gc.15}
 \end{equation}
 for all bounded measurable $f$.
  Hence $U_{\ep c,G}$ has a density  
   \bea
 &&
\hspace{-.3 in}u_{\ep c,G} (x,y)  =\sum_{n=0}^{\ff}\ep^{n }\,\,\int v_{\ep c}(x,z_{1})G(z_{1},\,dz_{2})v_{\ep c}( z_{2}, z_{3})\cdots \label{gc.3}\\
&&\hspace{-.1 in}\cdots G(z_{2n-3},\,dz_{2n-2})v_{\ep c}( z_{2n-2}, z_{2n-1})G(z_{2n-1},\,dz_{2n})v_{\ep c}( z_{2n}, y)\nn\\
&&\hspace{2.8 in}\prod_{j=1}^{n}c(z_{2j-1})\,dm(z_{2j-1}),\nn
 \eea
 with
    \begin{equation}
\sup_{z}\int   u_{\ep c,G}  ( z , y) dm(y) <\ff.\label{gc.4c}
 \end{equation}
  We can also write this as 
    \bea
 &&
\hspace{-.3 in}u_{\ep c,G} (x,y)   =\sum_{n=0}^{\ff}\ep^{n }\,\,\int v_{\ep c}(x,z_{1})j(z_{1},z_{2})v_{\ep c}( z_{2}, z_{3})\cdots \label{gc.3f}\\
&&\hspace{-.1 in}\cdots j(z_{2n-3},z_{2n-2})v_{\ep c}( z_{2n-2}, z_{2n-1})j(z_{2n-1},z_{2n})v_{\ep c}( z_{2n}, y)\prod_{j=1}^{2n}\,dm(z_{j}).\nn
 \eea

 Similar expansions can be given for the semigroup which has therefore a density:
    \bea
 &&
\hspace{-.3 in}p_{\ep c,G,t} (x,y)  =\sum_{n=0}^{\ff}\ep^{n }\,\,\int_{0\leq t_{1}\leq\cdots\leq t_{n}\leq t}  q_{t_1}(x,z_{1})c(z_1)G(z_{1},\,dz_{2})q_{t_2-t_1}( z_{2}, z_{3})\cdots\nn \\
&&\hspace{.5 in}\cdots c(z_{2n-1})G(z_{2n-1},\,dz_{2n})q_{t-t_n}( z_{2n}, y)\prod_{j=1}^{n}\,dm(z_{2j-1})\,dt_{j},\label{jgc.1}
 \eea
 where $q_t$ denotes the kernel of the semigroup associated with the process killed at rate $\ep c$.

 To see this let
 \begin{equation}
 W_{c,G,n,t}f(x) =: P^{x}\( f\(Y_{t}\);\,\tau(T_{n})<t<\tau(T_{n+1})\). \label{jgc.10a}
 \end{equation}
 Using the fact that
 $\tau_{t+u}=\tau_{t}+\tau_{u}\circ \th_{\tau_{t}}$ we have 
  \begin{eqnarray}
  && =P^{x}\( f\(Y_{t}\);\,\tau(T_{n})<t<\tau(T_{n+1})\Big | \tau(T_{n})\) 
 \label{jgc.10}\\ 
 &&=P^{x}\(\(f\(Y_{t-\tau (T_{n})}\);\,t-\tau (T_{n})<\tau(\la_{n+1})\)\circ \th_{\tau (T_{n})}\Big | \tau(T_{n})\)   \nonumber\\
&&=P^{x}\(P^{Y_{\tau(T_{n})}}\(f\(X_{t-\tau (T_{n})}\);\,t-\tau (T_{n})<\tau(\la_{n+1})\) \Big | \tau(T_{n})\)   \nonumber\\
&&=P^{x}\(\int G(Y^{-}_{\tau(T_{n})},\,dz )\right.\nn\\
&&\left.\hspace{.5 in}P_{\la_{n+1}}^{z }\(f\(X_{t-\tau (T_{n})}\);\,t-\tau (T_{n})<\tau(\la_{n+1})\)\Big | \tau(T_{n})\) .   \nonumber
 \end{eqnarray}
Conditional on $\tau(T_{n})$ we have 
\begin{eqnarray}
&&P_{\la}^{x}\( f\(X_{t-\tau (T_{n})}\);\,t-\tau (T_{n})<\tau(\la )\)
\label{j5.3}\\
&& =P_{\la}^{x}\( f\(X_{t-\tau (T_{n})}\);\,A_{t-\tau (T_{n})}< \la  \)\nonumber\\
&& =P^{x}\(  e^{-A_{t-\tau (T_{n})}}f\(X_{t-\tau (T_{n})}\)  \).  \nonumber
\end{eqnarray}
Hence setting  
\begin{equation}
P_{c,t}f(x)=P^{x}\(  e^{-A_{t}}f\(X_{t}\) \),\label{j5.6}
\end{equation}
 we have shown that
  \begin{eqnarray}
W_{c,G,n,t}f(x) &=& P^{x}\(\int G(Y^{-}_{\tau(T_{n})},\,dz )P_{c,t-\tau (T_{n}) }f (z )   \) 
 \label{jgc.11}\\ 
 &=&P^{x}\(  GP_{c,t-\tau (T_{n}) }f (Y^{-}_{\tau(T_{n})})   \).   \nonumber
 \end{eqnarray}
 
Using once again the fact that
 $\tau_{t+u}=\tau_{t}+\tau_{u}\circ \th_{\tau_{t}}$ and the strong  Markov property, we see that for any $h$
 \begin{eqnarray}
 &&
 P^{x}\(  h(Y^{-}_{\tau(T_{n})},t-\tau (T_{n}) ) \Big | \tau(T_{n-1})\)  \label{jgc.12}\\
 &&=P^{x}\(  \(h(Y^{-}_{\tau(\la_{n})},t-\tau (T_{n-1})-\tau(\la_{n})) \)\circ \th_{\tau (T_{n-1})}  \Big | \tau(T_{n-1})\) 
\nn\\
 &&= P^{x}\( P^{Y_{\tau(T_{n-1})}} \(h(Y^{-}_{\tau(\la_{n})},t-\tau (T_{n-1})-\tau(\la_{n})) \) \Big | \tau(T_{n-1})\)  \nonumber\\
 &&= P^{x}\( \int G(Y^{-}_{\tau(T_{n-1})},\,dz ) P^{z } \(h(X^{-}_{\tau(\la_{n})},t-\tau (T_{n-1})-\tau(\la_{n})) \)  \Big | \tau(T_{n-1})\) .  \nonumber
 \end{eqnarray}
  Using the change of variables formula, \cite[Chapter 6, (55.1)]{DM2}, and the fact that $X^{-}_{t}=X_{t}$ for a.e. t, we see that   conditionally on $\tau (T_{n-1})$ 
\bea
&&
P_{\la_{n}}^{z}\(h(X^{-}_{\tau_{\la_{n}}},t-\tau (T_{n-1})-\tau(\la_{n})) \)\label{j6.82adn}\\
&&=P_{\la_{n}}^{z}\(h(X_{\tau_{\la_{n}}},t-\tau (T_{n-1})-\tau(\la_{n})) \)\nn\\
&&=E^{x}\(\int_{0}^{\ff}e^{-r} h\(X_{\tau_{r}},t-\tau (T_{n-1})-\tau_{r}\)\,dr\)\nn\\
&&=
P^{z}\(\int_{0}^{A_{\ze}}e^{-r} h\(X_{\tau_{r}},t-\tau (T_{n-1})-\tau_{r}\)\,dr\)\nn\\
&&=P^{z}\(\int_{0}^{\ff}e^{-A_{s}} h\(X_{s},t-\tau (T_{n-1})-s\)\,dA_{s}\)\nn\\
&&=\int_{0}^{\ff} P^{z}\(e^{-A_{s}} h\(X_{s},t-\tau (T_{n-1})-s\)c(X_{s})\)\,ds\nn\\
&&=\int_{0}^{\ff}P_{c,s}(ch(\cdot,t-s-\tau (T_{n-1}) ))(z)\,ds. \nn
\eea
Combining this with (\ref{jgc.12}) we have 
 \bea
 &&
 P^{x}\(  h(Y^{-}_{\tau(T_{n})},t-\tau (T_{n}))   \)\label{jgc.13}\\
 &&=\int P^{x}\(  GP_{c,s_{n}}ch(Y^{-}_{\tau(T_{n-1})},t-s_{n}-\tau (T_{n-1}) )   \)\,ds_{n}.\nn 
 \eea
Iterating this we obtain,  conditionally on $\tau (T_{n-1})$, then $\tau (T_{n-2}), \cdots $
\begin{eqnarray}
&&
 P^{x}\(  h(Y^{-}_{\tau(T_{n})}),t-\tau (T_{n}))    \)\label{js1}\\
 &&=\int P^{x}\( GP_{c,s_{n}}ch(Y^{-}_{\tau(T_{n-1})},t-s_{n}-\tau (T_{n-1}) )   \)\,ds_{n}
\nn\\
 &&=\int P^{x}\(GP_{c,s_{n-1}}c   GP_{c,s_{n}}ch(Y^{-}_{\tau(T_{n-2})},t-s_{n}-s_{n-1}-\tau (T_{n-2}) )   \)\,ds_{n-1}\,ds_{n} \nonumber\\
 &&=\cdots  \nonumber\\
 &&=\int P^{x}\(GP_{c,s_{2}}c  \cdots GP_{c,s_{n}}ch(Y^{-}_{\tau(\la)},t-\sum_{j=2}^{n}s_{j}-\tau (\la) )   \)\prod_{j=2}^{n}\,ds_{j} \nonumber\\
 &&=\int  P_{c,s_{1}}c GP_{c,s_{2}}c  \cdots GP_{c,s_{n}}ch(x,t-\sum_{j=1}^{n}s_{j}  )    \prod_{j=1}^{n}\,ds_{j} \nonumber
\end{eqnarray}
where the last step used (\ref{j6.82adn}).
 Applying this to (\ref{jgc.11})  we have that
 \begin{equation}
W_{c,G,n,t}f(x)=\int  P_{c,s_{1}}c GP_{c,s_{2}}c  \cdots GP_{c,s_{n}}c GP_{c,t-\sum_{j=1}^{n}s_{j} } f(x )    \prod_{j=1}^{n}\,ds_{j}.\label{jgc.14}
 \end{equation}
(\ref{jgc.1}) then follows as in the proof of (\ref{gc.3}).

 Let
\begin{equation}
\hat{G}(x,dy)=\frac{1}{c(x)}j(y,x)m(dy).\label{}
\end{equation}
By  (\ref{jl.1}), $\hat{G}(x,dy)$ is a probability kernel on $S\times \BB (S)$. We now add jumps to the dual process $\wh{X}$, where  $c(x)$ will again govern the rate of jumps but we use   $\hat{G}(x,dy)$ to  describe the distribution of the jumps from position $x$. 
  Performing the same calculation as before, but with the dual process and using (\ref{gc.3f}) and  (\ref{swi}), we see that the two processes obtained by adding jumps have dual potential kernels:
  \begin{equation}
 \hat{u}_{\ep c,\hat{G}}  (x , y)=u_{\ep c,G}  ( y , x) \label{dual5}
  \end{equation}
  The same will be true for the associated resolvents and semigroups. The duality assumptions are verified and we can therefore define a loop measure associated with these processes.
 
  We use $\mu_{(\ep)}$ to denote the loop measure associated to $Y$, where we have replaced $c$ by $\ep c$.

 The next Lemma is needed for the proof of Theorem \ref{theo-jump}.
 
  \begin{lemma}\label{lem-norm}
  Assume (\ref{gc.4})  (which implies  (\ref{gc.4f})).
 Then for any positive measure $\nu$ and  $\ep$ sufficiently small.
 \begin{equation}
 \sup_{z}\int   u_{\ep c,G} ( z , y) d\nu(y)\leq C\sup_{z}\int   u ( z , y) d\nu(y),\label{gc.4r}
 \end{equation}
and
 \begin{equation}
\|\nu\|_{u_{\ep c,G}^2,\ff}\leq C\|\nu\|_{u^2,\ff},\label{gc.54}
\end{equation}
  
 \end{lemma}

  {\bf  Proof of Lemma \ref{lem-norm}: }
  (\ref{gc.4r}) follows immediately from (\ref{gc.3}) and (\ref{pd.1}).
  
It also follows from (\ref{gc.3}) that for a positive measure $\nu $
\bea
&&
\int  u^{2}_{\ep c,G} (x,y)\,d\nu(y)\label{gc.50} \\
&&=\sum_{m,n=0}^{\ff}\ep^{m+n } \int V_{\ep c}\(c GV_{\ep c}\)^{m}(x,y) V_{\ep c}\(c GV_{\ep c}\)^{n}(x,y)\,d\nu(y) \nn \\
&&=\sum_{m,n=0}^{\ff}\ep^{m+n } \int V_{\ep c}\(c GV_{\ep c}\)^{m-1}c G(x,\,dz_{1}) V_{\ep c}\(c GV_{\ep c}\)^{n-1}c G(x,\,dz_{2})\nn\\
&&\hspace{2 in}\(\int v_{\ep c}( z_{1},y)v_{\ep c}( z_{2},y)\,d\nu(y)\). \nn
\eea
Hence for $\ep$ small enough
\bea
&&
\sup_{x}\int  u^{2}_{\ep c,G} (x,y)\,d\nu(y)\label{gc.51} \\
&&\leq \sum_{m,n=0}^{\ff}\ep^{m+n } \sup_{x}\int V_{\ep c}\(c GV_{\ep c}\)^{m-1}c G(x,\,dz_{1}) V_{\ep c}\(c GV_{\ep c}\)^{n-1}c G(x,\,dz_{2})\nn\\
&&\hspace{2 in}\sup_{z_{1},z_{2}}\(\int v_{\ep c}( z_{1},y)v_{\ep c}( z_{2},y)\,d\nu(y)\) \nn\\
&&\leq C\|\nu\|_{u^2,\ff}. \nn
\eea
 Similarly
 \bea
&&
\int  u^{2}_{\ep c,G} (y,x)\,d\nu(y)\label{gc.52} \\
&&=\sum_{m,n=0}^{\ff}\ep^{m+n } \int V_{\ep c}\(c GV_{\ep c}\)^{m}(y,x) V_{\ep c}\(c GV_{\ep c}\)^{n}(y,x)\,d\nu(y) \nn \\
&&=\sum_{m,n=0}^{\ff}\ep^{m+n } \int \(c GV_{\ep c}\)^{m}(z_{1},x)\(c GV_{\ep c}\)^{n}(z_{2},x)
\,dm(z_{1})\,dm(z_{2})\nn\\
&&\hspace{2 in}\(\int v_{\ep c}(  y,z_{1})v_{\ep c}( y,z_{2})\,d\nu(y)\). \nn
\eea
Using (\ref{gc.4f}) it follows that for $\ep$ small enough
\bea
&&
\sup_{x}\int  u^{2}_{\ep c,G} (x,y)\,d\nu(y)\label{gc.53} \\
&&\leq \sum_{m,n=0}^{\ff}\ep^{m+n } \sup_{x} \int \(c GV_{\ep c}\)^{m}(z_{1},x)\(c GV_{\ep c}\)^{n}(z_{2},x)
\,dm(z_{1})\,dm(z_{2})\nn\\
&&\hspace{2 in}\sup_{z_{1},z_{2}}\(\int v_{\ep c}(  y,z_{1})v_{\ep c}( y,z_{2})\,d\nu(y)\) \nn\\
&&\leq C\|\nu\|_{u^2,\ff}. \nn
\eea
\qed

It follows from  \cite[Lemma 3.3]{LMR}   that $\|\nu\|_{u^2,\ff}$
 is a proper norm for $u_{\ep c,G, n}$.

 Theorem \ref{theo-jump}  follows from the next Lemma.
 
 \begin{lemma}\label{lem-jump} Under the assumptions of Lemma \ref{lem-norm}, for any  $ \nu_{1},\ldots, \nu_{k}\in \mathcal{R}_{\|\cdot\|_{u^2,\ff}}^{+} $
\begin{eqnarray}
&&{d \over d\ep}\int \prod_{j=1}^{k}u_{\ep c,G}(y_{j},y_{j+1})  \prod_{j=1}^{k} \,d\nu_{j}(y_{j})\,\,\,\Big|_{\ep=0}
\label{f.1}\\
&&   =\sum_{i=1}^{k}\int\(\prod_{j=1}^{i-1}u (y_{j},y_{j+1})\)\(- \int u(y_{i},x)u(x,y_{i+1})c(x)\,   dm(x)\right.\nn\\
&&\left.\hspace{1.4 in}  +\int_{S\times S}u(y_{i},z_{1})c(z_{1}) G(z_{1},\,dz_{2}) u(z_{2},y_{i+1})\,dm (z_{1})\)\nn\\
&& \hspace{2.5 in}    \(  \prod_{j=i+1}^{k}u (y_{j},y_{j+1})\)  \prod_{j=1}^{k} \,d\nu_{j}(y_{j}),  \nonumber
\end{eqnarray}
with $y_{k+1}=y_{1}$.
\end{lemma}

 {\bf  Proof of Lemma \ref{lem-jump}: }
 Set
 \begin{equation}
 \mathcal{I}(\ep)=\int \prod_{j=1}^{k}u_{\ep c,G}(y_{j},y_{j+1})  \prod_{j=1}^{k} \,d\nu_{j}(y_{j})\label{f.2}
 \end{equation}
 with $y_{k+1}=y_{1}$.

We can write (\ref{gc.3}) as
 \bea
 &&
\hspace{-.3 in} u_{\ep c,G} (x,y) \label{gc.5}\\
 &&\hspace{-.3 in}=v_{\ep c}(x,y)+\ep \int v_{\ep c}(x,z_{1})c(z_{1})G(z_{1},\,dz_{2})v_{\ep c}( z_{2}, y)\,dm(y_{1})\nn\\
 &&\hspace{-.3 in}+\ep ^{2}\int v_{\ep c}(x,z_{1})c(z_{1})G(z_{1},\,dz_{2})v_{\ep c}( z_{2}, y)c(z_{2})G(z_{2},\,dz_{3}) u_{\ep c,G}( z_{3}, y)\,dm(y_{1})\,dm(y_{3})\nn\\
 &&\hspace{-.3 in}=v_{\ep c}(x,y)+\ep V_{\ep c} c GV_{\ep c}   (x,y)+\ep ^{2}V_{\ep c}c GV_{\ep c} c GU_{\ep c,G} (x,y),\nn
 \eea
 with operator notation.
 
 We substitute this in (\ref{f.2}) and collect terms to obtain
  \begin{equation}
 \mathcal{I}(\ep)=I(\ep)+\ep \sum_{i=1}^{k}\mathcal{J}_{i}(\ep)+\mathcal{K}(\ep),\label{f.3}
 \end{equation}
where
\begin{equation}
I(\ep)=\int \prod_{j=1}^{k}v_{\ep c}(y_{j},y_{j+1})  \prod_{j=1}^{k} \,d\nu_{j}(y_{j})\label{f.4}
\end{equation}
\begin{eqnarray}
&&\hspace{-.3 in}\mathcal{J}_{i}(\ep)
\label{f.5}\\
&&\hspace{-.3 in} =\int \( \prod_{j=1}^{i-1}v_{\ep c}(y_{j},y_{j+1})\)V_{\ep c} c GV_{\ep c}(y_{i},y_{i+1})   \( \prod_{j=i+1}^{k}v_{\ep c}(y_{j},y_{j+1})\)  \prod_{j=1}^{k} \,d\nu_{j}(y_{j}),  \nonumber
\end{eqnarray}
and $\mathcal{K}(\ep)$ represents all the remaining terms. Noting that $ \mathcal{I}(0)=I(0)$, we can write (\ref{f.3}) as
 \begin{equation}
 \mathcal{I}(\ep)-\mathcal{I}(0)=I(\ep)-I(0)+\ep \sum_{i=1}^{k}\mathcal{J}_{i}(\ep)+\mathcal{K}(\ep). \label{f.6}
\end{equation}

Note also that $I(\ep)$ of (\ref{f.4}) is a special case of  the $I(\ep)$ of (\ref{rg.1}) with $\nu (dx)=c(x)\,dm(x)$. Hence by (\ref{df.7})
\begin{eqnarray}
&&\lim_{\ep\to 0}{I(\ep)-I(0) \over \ep}
\label{f.7}\\
&&=  -\sum_{i=1}^{k}\int\(\prod_{j=1}^{i-1}u (y_{j},y_{j+1})\)\( \int u(y_{i},x)u(x,y_{i+1})c(x)\,dm(x)\)\nn\\
&& \hspace{2.3 in}    \(  \prod_{j=i+1}^{k}u (y_{j},y_{j+1})\)  \prod_{j=1}^{k} \,d\nu_{j}(y_{j}).  \nonumber
\end{eqnarray}

Since $v_{\ep c}(x,y )\uparrow u(x,y )$ as $\ep\downarrow 0$, it follows by the Monotone Convergence Theorem that
\begin{eqnarray}
&&\lim_{\ep\to 0} \mathcal{J}_{i}(\ep)  
\label{f.75}\\
&&=  \int\(\prod_{j=1}^{i-1}u (y_{j},y_{j+1})\) \nn\\
&& \hspace{1.4 in}  \(\int_{S\times S}u(y_{i},z_{1})c(z_{1}) G(z_{1},\,dz_{2}) u(z_{2},y_{i+1})\,dm (z_{1})\)\nn\\
&& \hspace{2.3 in}    \(  \prod_{j=i+1}^{k}u (y_{j},y_{j+1})\)  \prod_{j=1}^{k} \,d\nu_{j}(y_{j}).  \nonumber
\end{eqnarray}

To complete the proof of our Lemma it remains to show that
\begin{equation}
\lim_{\ep\to 0}{\mathcal{K}(\ep) \over \ep}=0.\label{f.8}
\end{equation}
 However every term in $\mathcal{K}(\ep)$ comes with a pre-factor of $\ep^{m}$ for some $m\geq 2$, so we need only bound the integrals uniformly in $\ep$. For this we will use Lemma \ref{lem-norm}. We illustrate this with the most complicated term, which has the  pre-factor   $\ep^{2k}$:
\begin{eqnarray}
&&\int \prod_{j=1}^{k}V_{\ep c}c GV_{\ep c} c GU_{\ep c,G} (y_{j},y_{j+1})  \prod_{j=1}^{k} \,d\nu_{j}(y_{j})
\label{f.9}\\
&&=  \int v_{\ep c}(y_{1},x_{1})c GV_{\ep c} c GU_{\ep c,G}(x_{1},y_{2}) \cdots \nn\\
 &&\hspace{.3 in} v_{\ep c}(y_{k-1},x_{k-1})c GV_{\ep c} c GU_{\ep c,G}(x_{ k-1},y_{ k})\, \,\nn\\
 &&\hspace{1.4 in}v_{\ep c}(y_{k},x_{k})c GV_{\ep c} c GU_{\ep c,G}(x_{k},y_{1})  \prod_{j=1}^{k} \,d\nu_{j}(y_{j}) \,dm (x_{j})\nn\\
  &&=  \int c GV_{\ep c} c GU_{\ep c,G}(x_{1},y_{2}) v_{\ep c}(y_{2},x_{2})\cdots \nn\\
 &&\hspace{.3 in} c GV_{\ep c} c GU_{\ep c,G}(x_{ k-1},y_{ k})\,  v_{\ep c}(y_{k},x_{k})\,\nn\\
 &&\hspace{1.3 in}c GV_{\ep c} c GU_{\ep c,G}(x_{k},y_{1})  v_{\ep c}(y_{1},x_{1})\prod_{j=1}^{k} \,d\nu_{j}(y_{j}) \,dm (x_{j}),\nn
 \eea
 where the last step is just a rearrangement. We can rewrite this as 
 \begin{eqnarray}
 && \int \(\int c GV_{\ep c} c GU_{\ep c,G}(x_{1},y_{2})v_{\ep c}(y_{2},x_{2})\,d\nu_{2}(y_{2})\)\cdots  \label{f.10}\\
 &&\hspace{.3 in}\( \int c GV_{\ep c} c GU_{\ep c,G}(x_{ k-1},y_{ k})\, v_{\ep c}(y_{k},x_{k})\,d\nu_{k}(y_{k})\)\,\nn\\
 &&\hspace{1 in}\(\int c GV_{\ep c} c GU_{\ep c,G}(x_{k},y_{1})v_{\ep c}(y_{1},x_{1})\,d \nu_{1}(y_{1}) \)\,\,\prod_{j=1}^{k} \,\,dm (x_{j}). \nn 
 \end{eqnarray} 
 Then by Lemma \ref{lem-norm}, and the fact that $v_{\ep c}\leq u$ 
 \begin{eqnarray}
 &&\hspace{-.3 in}\int c GV_{\ep c} c GU_{\ep c,G}(x_{1},y_{2})v_{\ep c}(y_{2},x_{2})\,d\nu_{2}(y_{2})
 \label{f.11}\\
 &&\hspace{-.3 in}=\int c(x_{1}) G(x_{1}, dz_{1})v_{\ep c}(z_{1},z_{2}) c(z_{2} )G(z_{2}, dz_{3})u_{\ep c,G}(z_{3},y_{2})v_{\ep c}(y_{2},x_{2})\,d\nu_{2}(y_{2}) \,dm (z_{2})  \nonumber\\
 &&\hspace{-.3 in}\leq \int c(x_{1}) G(x_{1}, dz_{1})v_{\ep c}(z_{1},z_{2}) c(z_{2} )G(z_{2}, dz_{3})\,dm (z_{2})  \nonumber\\
 &&\hspace{2 in}\sup_{z_{3},x_{2}}\int u_{\ep c,G}(z_{3},y_{2})v_{\ep c}(y_{2},x_{2})\,d\nu_{2}(y_{2}) \nn\\
 &&\hspace{-.3 in}\leq C \|\nu_{2}\|_{u^2,\ff} \int c(x_{1}) G(x_{1}, dz_{1})v_{\ep c}(z_{1},z_{2}) c(z_{2} )G(z_{2}, dz_{3})\,dm (z_{2}).  \nonumber\\
 \end{eqnarray}
Using (\ref{gc.4b}), the fact that $c$ is bounded  and the fact that $G(\cdot, dz)$ is a probability density 
 \begin{equation}
  \int c(x_{1}) G(x_{1}, dz_{1})v_{\ep c}(z_{1},z_{2}) c(z_{2} )G(z_{2}, dz_{3})\,dm (z_{2})\leq C c(x_{1}).\label{f.12}
 \end{equation}
 Thus (\ref{f.10}) is bounded  independently of $\ep$ by
 \begin{equation}
C\int  \prod_{j=1}^{k}c(x_{j}) \,\,dm (x_{j})<\ff \label{f.13}
 \end{equation}
 since  $c(x)$ is integrable.
 
 The other terms of $\mathcal{K}(\ep)$ can be bounded similarly.\qed

\def\noopsort#1{} \def\printfirst#1#2{#1}
\def\singleletter#1{#1}
            \def\switchargs#1#2{#2#1}
\def\bibsameauth{\leavevmode\vrule height .1ex
            depth 0pt width 2.3em\relax\,}
\makeatletter
\renewcommand{\@biblabel}[1]{\hfill#1.}\makeatother
\newcommand{\bysame}{\leavevmode\hbox to3em{\hrulefill}\,}

\def\wh{\widehat}
\def\ol{\overline}

{\footnotesize

\noindent
\begin{tabular}{lll} &  Yves Le Jan & \hskip20pt  Jay Rosen \\   &  Equipe Probabilit\'es et Statistiques
     & \hskip20pt Department of Mathematics \\    &
  Universit\'e Paris-Sud, B\^atiment 425
     & \hskip20pt  College of Staten Island, CUNY  \\ &  91405 Orsay Cedex
     & \hskip20pt  Staten Island, NY
10314
 \\    &   France
     & \hskip20pt U.S.A. \\   &   
yves.lejan@math.u-psud.fr 
     & \hskip20pt    jrosen30@optimum.net
\end{tabular}
\bigskip

}

  \end{document}